\documentclass[letterpaper,11pt]{article}

\usepackage{amsmath,amsfonts,amsthm,amstext,amssymb,amscd,amsbsy}
\usepackage[latin1]{inputenc}
\usepackage[T1]{fontenc}
\usepackage{indentfirst}
\usepackage[dvips]{color}
\usepackage{epsfig}
\usepackage{color}
\usepackage{graphics}
\usepackage[normalem]{ulem}

\newtheorem{defi}{Definition}[section]
\newtheorem{cor}[defi]{Corollary}
\newtheorem{exem}[defi]{Example}
\newtheorem{prop}[defi]{Proposition}
\newtheorem{obs}[defi]{Remark}
\newtheorem{lema}[defi]{Lemma}
\newtheorem{teo}[defi]{Theorem}

\newcommand{\R}{\mathbb{R}}
\newcommand{\C}{\mathbb{C}}

\newcommand{\F}{\mathbb{F}}
\newcommand{\fimdem}{\hfill\rule{0.25cm}{0.25cm}}
\newcommand{\seta}{\longrightarrow}
\newcommand{\Ad}{\mbox{Ad}}
\newcommand{\ad}{\mbox{ad}}
\newcommand{\opbra}{\left\{}
\newcommand{\clbra}{\right\}}
\newcommand{\gotico}{\mathfrak}

\newcommand{\dem}{\noindent \textbf{Proof:}\quad}

\newcommand{\ep}{\varepsilon}
\def\bq{\begin{equation}}
\def\eq{\end{equation}}

\title{Equigeodesics on Flag Manifolds}
 \author{Nir Cohen
\footnote{Department of Applied Mathematics; corresponding author,
nir@ime.unicamp.br},
 Lino Grama
\footnote{Department of Mathematics, linograma@gmail.com}, and
 Caio J.C. Negreiros
\footnote{Department of Mathematics, caione@ime.unicamp.br},\\
Institute of Mathematics, Statistics and Scientific Computation,\\
P.O. Box    6065, University of Campinas - UNICAMP.}
\date{}

\begin{document}

\maketitle

%{ \footnotesize \parindent

\begin{abstract}
{\em This paper provides a characterization of homogeneous curves on
a geometric flag manifold which are geodesic with respect to any
invariant metric. We call such curves homogeneous equigeodesics. We
also characterize homogeneous equigeodesics whose associated Killing
field is closed, hence, the corresponding geodesics is closed.}
\end{abstract}



% \tableofcontents

\section{Introduction}

Let $(M,g)$ be a Riemannian manifold and let $\gamma$ be a geodesic
passing at $p\in M$ with direction vector $X\in T_pM.$ The geodesic
$\gamma$ is called {\em homogeneous} if it is the orbit of a
1-parameter subgroup of $G$, that is, $\gamma(t)=\exp_p(tX).$
Furthermore, $(M,g)$ is called a {\it g.o. manifold} (geodesic orbit
manifold) if {\em every} geodesic is homogeneous.

The g.o. property is particularly meaningful if we restrict the
discussion to homogeneous spaces $M=G/K$ and $G$-invariant metrics
$g$. In this case, we may choose $p$ to be the origin, i.e. the
trivial coset, and identify $T_pM$ with the corresponding subspace of the Lie algebra
$\gotico{g}.$ The set of g.o. manifolds includes all the symmetric
spaces; their classification up to dimension 6 can be found in
Kowalski and Vanhecke \cite{K2}.

The normal metric is g.o. on any flag manifold \cite{cheeger}.
Alekseevsky and Arvanitoyeorgos \cite{arva3} showed that the only
flag manifolds which admit a g.o. metric not homothetic to the
normal metric are $SO(2l+1)/U(l)$ and $Sp(l)/U(1)\times Sp(l-1)$.
More recently, Alekseevsky and Nikonorov  \cite{alekniko} obtained a
classification of compact, simply-connected homogeneous g.o. spaces
with positive Euler characteristic.

According to Kowalski and Szenthe \cite{K1}, every homogeneous
Riemannian manifold admits homogeneous geodesics. In the present
paper we show that every flag manifold of the $A_l$ type admits
homogeneous {\it equi}geodesics, namely homogeneous curves
$\gamma$ which are geodesic with respect to {\it any}
$G$-invariant metric. We shall give a full characterization of
homogeneous equigeodesics $\gamma$ in terms of the corresponding
vectors $X$, which we call {\em equigeodesic vectors}. Our
starting point is the following algebraic characterization:

\begin{teo}\label{first}
A tangent vector $X$ for the flag manifold $\F(n;n_1,\cdots,n_k)$
is equigeodesic iff $[X,\Lambda X]_\gotico{m}=0$  for every
invariant metric $\Lambda.$\end{teo}
\noindent By the classical adjoint representation, $X$ corresponds
to an $n\times n$ skew-Hermitian matrix $A$ with blocks $a_{ij}\in
M_{n_i,n_j}(\C)$, with $a_{ii}=0;$ similarly, the metric $g$ is
represented by a symmetric $n\times n$ matrix $\Lambda$ with
positive entries $\lambda_{ij}$, constant in each block, with
$\lambda_{ii}=0$. The inner product is $g(X,Y)=(\Lambda X,Y)$ where
the product $\Lambda X$ is the Hadamard (or termwise) product
(\cite{hj}). In these terms we show the following result. Recall
that the vector $X$ extends uniquely to a Killing field which
contains $\gamma$ as a trajectory.  If the Killing field is closed, by definition $\gamma$
is closed, but the converse need not hold. 

\begin{teo}\label{second}
\quad (i) $X$ is equigeodesic iff $a_{ij}a_{jm}=0$ for all $i,j,m$
distinct, $1\leq i,j,m \leq k$.\\
%
%\noindent {(ii) The associated equigeodesic $\gamma$ is closed iff
%the eigenvalues of $A$ are pairwise $\mathbb{Q}$-linearly
%dependent.}
%
(ii) The eigenvalues of $A$ are commensurate iff $X$ defines a
closed Killing field.
\end{teo}

\noindent We show item (ii) by putting the matrix $A$ in an
essentially diagonal canonical form, and then using a recent
characterization of closed Killing fields (Flores et al.,
\cite{Piccione}).

%$\gamma$ is closed iff {\em
%every} trajectory of the Killing field is closed [??] {\bf
%encontrar fonte!}

In the special case of the full flag manifold, where blocks of $A$
and $\Lambda$ are scalar, Theorem \ref{second} simplifies
considerably:

\begin{cor}\label{cor3}
(i) $X$ is equigeodesic in $\F(n)$ iff $A$ is permutation-similar to
a diagonal matrix.\quad
%
%\noindent \sout{(ii) The corresponding geodesic is closed iff these
%entries are pairwise $\mathbb{Q}$-linearly dependent.}
%
(ii) $\gamma$ is closed if the entries (rather than eigenvalues) of
$A$ are commensurate.\end{cor}

The simplest equigeodesic choice is $X\in\gotico{u}_\alpha$, where
$A$ has a single pair of non-zero entries. The resulting geodesic
$\gamma$ is closed, and a simpler argument suffices to prove its
closure. Indeed, $\gamma$ is embedded in a totally geodesic 2-sphere
$S^2$ embedded in $\F(n)$, having $\gotico{u}_\alpha$ as a tangent
space. Thus $\gamma=S^1,$ a closed geodesic.

In this construction, the curve $\gamma$ and surface $S^2$ are
both {\em equiharmonic} (for this notion, see Black \cite{black})
in $\F(n)$. We mention that a class of equiharmonic maps from
$S^2$ to $\F(n)$ was found by Negreiros in \cite{Neg1}; and it is
still an open problem whether any harmonic map between $S^2$ and
$\F(n)$ is necessarily equiharmonic. In this paper we have shown
that a homogeneous geodesic curve need not be equigeodesic. We are now studying equigeodesics in flag manifolds of other Lie groups (classical and exceptional). 

\section{The geometry of flag manifolds}

In this section we briefly review basic facts on the structure of
homogeneous spaces and flag manifolds; and describe the $T$-roots
system used in constructing the partial flag manifold
$\F(n;n_1,\cdots,n_k)$.\\

{\em I. Homogeneous spaces.} Consider the homogeneous manifold
$M=G/K$ with $G$ a compact semi-simple Lie group and $K$ a closed
subgroup. Let $\gotico{g}$ and $\gotico{k}$ be the corresponding
Lie algebras. The Cartan-Killing form $\left\langle
,\right\rangle$ is nondegenerate and negative definite in
$\gotico{g}$, thus giving rise to the direct sum decomposition
$\gotico{g=k\oplus m}$ where $\gotico{m}$ is $\Ad(K)$-invariant.
%that is, $\Ad(k)\gotico{m}\subset\gotico{m}$ for all $k\in K$.
We may identify $\gotico{m}$ with the tangent space $T_oM$ at
$o=eK$. The isotropy representation of a reductive homogeneous
space is the homomorphism $j:K\seta GL(T_oM)$ given by
$j(k)=\Ad(k)\big|_{\gotico{m}}$.

%There exists a one-to-one correspondence between $G$-invariant
%metrics $g$ on $M$ and $\Ad(K)$-invariant scalar products $B$ on
%$\gotico{m}$, that is,
% $$B(X,Y)=B(\Ad(k)X,\Ad(k)Y)\quad \mbox{\bf algo esta errado
%aqui???}$$ for all $X,Y\in \gotico{m}$ and $k\in K$, see
%\cite{KN}. {\bf nir: sugiro tirar}.

A metric $g$ on $M$ is defined by a scalar product on $\gotico{m}$
has the form $B(X,Y)=-\left\langle \Lambda X,Y\right\rangle$, with
$\Lambda:\gotico{m}\seta\gotico{m}$ positive definite with respect
to the Cartan-Killing form, see for example \cite{cheeger}. We
denote by $ds^2_{\Lambda}$ the invariant metric given by
$\Lambda$. We abuse of notation and say that $\Lambda$ itself is
an invariant metric.\\

{\em II. Generalized flag manifolds.} A homogeneous space $F=G/K$
is called a generalized flag manifold if $G$ is simple and the
isotropy group $K$ is the centralizer of a one-parametrer subgroup
of $G$, $\exp tw$ $(w\in \gotico{g})$. Equivalently, $F$ is an
adjoint orbit $\Ad(G)w$, where $w\in \gotico{g}$. The generalized
flag manifolds (also refered to as a K\"ahlerian $C$-spaces) have
been classified in \cite{BFR86},\cite{Wa}.

Here the direct sum decomposition
$\gotico{g}=\gotico{k}\oplus\gotico{m}$ has a more complete
description (see e.g. \cite{ale1},\cite{arva1}). Let
$\gotico{h}^{\C}$ be a Cartan subalgebra of the complexification
$\gotico{k}^{\C}$ of $\gotico{k}$, which is also a Cartan
subalgebra of $\gotico{g}^{\C}$. Let $R$ and $R_K$ be the root
systems of $\gotico{g}^{\C}$ and $\gotico{k}^{\C}$, respectively,
and $R_M=R\backslash R_K$ be the set of complementary roots. We
have the Cartan decompositions $$\gotico{g}^{\C}=\gotico{h}^{\C}
\oplus \sum_{\alpha\in R}{\gotico{g}_{\alpha}},\hspace{1cm}
\gotico{k}^{\C}=\gotico{h}^{\C} \oplus \sum_{\alpha\in
R_K}{\gotico{g}_{\alpha}}, \hspace{1cm}
\gotico{m}^{\C}=\sum_{\alpha\in R_M}{\gotico{g}_{\alpha}}$$
where $\gotico{m}^{\C}$ is isomorphic to $(T_oF)^{\C}$ and
$\gotico{h}=\gotico{h}^{\C}\cap\gotico{g}$. Thus, the real tangent
space of $T_oF$ is naturally identified with
 $$\displaystyle{\gotico{m}=\bigoplus_{\alpha\in R_M^+}
 \gotico{u}_{\alpha}}.$$
Unless $F$ is a full flag manifold, some of the spaces
$\gotico{u}_\alpha$ are not $\Ad(K)$-modules. To get the {\em
irreducible} $Ad(K)$-modules, we proceed as in \cite{ale1} or
\cite{arva2}. Let
$$\gotico{t}=Z(\gotico{k}^{\C})\cap\gotico{h}=\opbra
X\in\gotico{h}:\phi(x)=0 \,\, \forall \phi \in R_K\clbra.$$ If
$\gotico{h}^*$ and $\gotico{t}^*$ are the dual space of
$\gotico{h}$ and $\gotico{t}$ respectively, we consider the
restriction map
\begin{equation}
\label{projkappa} \kappa:\gotico{h}^*\seta \gotico{t}^*,
\quad\quad\quad\quad \kappa(\alpha)=\alpha|_{\gotico{t}}
\end{equation}
and set $R_T=\kappa(R_M)$. This set satisfies the axioms of a not
necessarily reduced root system, and its elements are called
$T$-roots. The irreducible $\ad(\gotico{k}^{\C})$-invariant
sub-modules of $\gotico{m}^{\C},$ and the corresponding
irreducible sub-modules for the $\ad(\gotico{k})$-module
$\gotico{m},$ are given by
 $$\gotico{m}^\C_{\xi}=\sum_{\kappa(\alpha)=\xi}
 {\gotico{g}_{\alpha}}\quad\quad(\xi\in
 R_T),\quad\quad\quad\quad
\gotico{m}_{\eta}=\displaystyle{\sum_{\kappa(\alpha)=\eta}{\gotico{u}_{\alpha}}}
 \quad\quad(\eta\in R^+_T).$$
Hence we have the direct sum of complex and real irreducible
modules,
  $$\gotico{m}^\C=\displaystyle{\sum_{\eta\in
R_T}{\gotico{m}^\C_{\eta}}},\quad\quad\quad
\gotico{m}=\displaystyle{\sum_{\eta\in
R^+_T}{\gotico{m}_{\eta}}}.$$

We fix a Weyl basis in $\gotico{m}^\C$, namely, elements
$X_{\alpha}\in \gotico{g}_{\alpha}$ such that $\left\langle
X_{\alpha},X_{-\alpha}\right\rangle=1$ and
$[X_{\alpha},X_{\beta}]=m_{\alpha,\beta}X_{\alpha+\beta}$, with
$m_{\alpha,\beta}\in \R$, $m_{\alpha,\beta}=-m_{\beta,\alpha}$,
$m_{\alpha,\beta}=-m_{-\alpha,-\beta}$ and $m_{\alpha,\beta}=0$ if
$\alpha+\beta$ is not a root. The corresponding {\em real} Weyl
basis in $\gotico{m}$ consists of the vectors
$A_{\alpha}=X_{\alpha}-X_{-\alpha}$,
$S_{\alpha}=i(X_{\alpha}+X_{-\alpha})$ and
$\gotico{u}_{\alpha}=\mbox{span}_{\R}\opbra
A_{\alpha},S_{\alpha}\clbra$, where $\alpha\in R^+$, the set of
positive roots.\\

An invariant metric $g$ on $F$ is uniquely defined by a scalar
product $B$ on $\gotico{m}^\C$ of the form
$$B(\cdot\,,\cdot)=-\left\langle \Lambda \cdot\,,\cdot
\right\rangle= \lambda_1(-\left\langle \cdot\,,\cdot
\right\rangle)|_{\gotico{m}_1}+\ldots + \lambda_j(-\left\langle
\cdot\,,\cdot \right\rangle)|_{\gotico{m}_j,}$$ where
$\lambda_i>0$ and $\gotico{m}_i$ are the irreducible
$\Ad(K)$-sub-modules. Each $\gotico{m}_i$ is an eigenspace of
$\Lambda$ corresponding to the eigenvalue $\lambda_i$. In
particular, the vectors $A_{\alpha},S_{\alpha}$ of the real Weyl
basis are eigenvectors of $\Lambda$ correponding to the same
eigenvalue $\lambda_{\alpha}$.\\

{\em III. Generalized flag manifolds of the geometric (or $A_l)$
type.} These are the spaces of type
$$\F(n;n_1,\ldots,n_s)={SU(n)}\; / \; {S(U(n_1)\times\ldots\times
U(n_s))},$$ where $n=n_1+\ldots +n_s$. Our description of
$T$-roots for these spaces follows \cite{arva2}.

The complexification of the real Lie algebra $\gotico{su}(n)$ is
$\gotico{sl}(n,\C)$. The Cartan sub-algebra of $\gotico{sl}(n,\C)$
can be identified with
$\gotico{h}=\opbra\mbox{diag}(\ep_1,\ldots,\ep_n);\ep_i\in\C,\sum{\ep_i}=0\clbra$.
The root system of the Lie algebra of $\gotico{sl}(n)$ has the
form $R=\opbra\alpha_{ij}=\ep_i-\ep_j:i\neq j\clbra$ and the
subset of positive roots is  $R^+=\opbra \alpha_{ij}:i<j\clbra$.
We have
 $$\begin{array}{cl} R_K=\{ \ep^i_a-\ep^i_b: &
 1\leq a \neq b\leq n_i\},\\
 R_K^+=\{ \ep^i_a-\ep^i_b: &
 1\leq a < b \leq n_i\},\\
 R_M^+=\{ \ep^i_a-\ep^j_b: & i<j,1\leq a\leq
n_i,1\leq b \leq n_j \},\end{array}$$ where we use the notation
$\ep^i_a=\ep_{n_1+\ldots+n_{i-1}+a}$. The sub-algebra $\gotico{t}$
of $\gotico{h}$ used in the construction of T-roots consists of
positive diagonal matrices of the form $diag\{\lambda_i
I_{n_i}\}_{i=1}^s.$ We conclude that the number of irreducible
$\Ad(K)$-submodules of $\F(n;n_1+\ldots+n_s)$ is
$\frac{1}{2}s(s-1)$. In the special case of the full flag manifold
$\F(n):=\F(n;1,\cdots,1)$, the sets of roots and $T$-roots
coincide.

\section{Equigeodesics on flag manifolds}

With these preliminaries we can now discuss in full detail the
characterization of equigeodesic vectors.

\begin{defi}
\label{defihomo} Let $(M=G/K,g)$ be a homogeneous Riemannian
manifold. A geodesic $\gamma(t)$ on $M$ through the origin $o$ is
called \textit{homogeneous} if it is the orbit of a $1$-parameter
subgroup of $G,$ that is, $$ \gamma(t)=(\exp tX)\cdot o, $$ where
$X \in \gotico{g}$. The vector $X$ is called a geodesic vector.
\end{defi}

Definition \ref{defihomo} establishes a 1:1 correspondence between
geodesic vectors $X$ and homogeneous geodesics at the origin. A
result of Kowalski and Vanhecke \cite{K2} implies, as a special
case, the following algebraic characterization.
\begin{teo}\label{hecke}If $g$ is a $G$-invariant metric, a vector $X \in
\gotico{g}\setminus\opbra 0 \clbra$ is a geodesic vector iff
\begin{equation}
\label{eqnfund} g(X_{\gotico{m}},[X,Z]_{\gotico{m}})=0,
\end{equation}
for all $Z\in \gotico{m}$.
\end{teo}
The following existence result is of interest:
\begin{teo}[\cite{K1}]
\label{givevector} If $G$ is semi-simple then $M$ admits at least
$m=dim(M)$ mutually orthogonal homogenous geodesics through the
origin $o$.
\end{teo}

An example is the classical flag manifold $\F(n)$ of real
dimension $n(n-1)$ and the real Weyl basis $\{
A_{\alpha},S_{\alpha},$ $\alpha \in R^+\}$ of the same size.
Actually, these vectors are geodesic vectors with respect to {\it any}
invariant metric $\Lambda$ on $\gotico{m}$, motivating the
following definition.

\begin{defi}
A curve $\gamma$ on $G/H$ is an {\normalfont equigeodesic} if it is a geodesic for
any invariant metric $ds^2_{\Lambda}$. If the equigeodesic is of the
form $\gamma(t)=(\exp tX)\cdot o$,  where $X \in \gotico{g}$, we say that
$\gamma$ is a {\normalfont homogeneous equigeodesic} and the vector $X$ is an
equigeodesic vector.
\end{defi}

Theorem \ref{hecke} simplifies in the special case of flag
manifolds and equigeodesic vectors.

\begin{prop}
\label{proplegal} Let $F$ be a flag manifold, with $\gotico{m}$
isomorphic to $ T_o F$. A vector $X\in \gotico{m}$ is equigeodesic
iff
\begin{equation}
\label{eqnlegal} [X,\Lambda X]_{\gotico{m}}=0,
\end{equation}
for any invariant metric $\Lambda$.
\end{prop}
\dem Let $g$ be the metric associated with $\Lambda$. For $X,Y\in
\gotico{m}$ we have
\begin{eqnarray*}
g(X,[X,Y]_{\gotico{m}}) &=& -\left\langle \Lambda
X,[X,Y]_{\gotico{m}} \right\rangle = - \left\langle \Lambda
X,[X,Y] \right\rangle = -\left\langle [X,\Lambda X],Y
\right\rangle,
\end{eqnarray*}
since the decomposition $\gotico{g=m+h}$ is $<,>$-orthogonal and the
Killing form is $\Ad(G)$-invariant, i.e., $\ad(X)$ is skew-Hermitian
with respect to $<,>$. Therefore $X$ is equigeodesic iff $[X,\Lambda
X]_{\gotico{m}}=0$ for any invariant scalar product $\Lambda$.
\fimdem
\\

In the ensuing analysis we assume $F=\F(n;n_1,\cdots,n_s),$ a
geometric flag manifold, and use the classical adjoint
representation of $G$ to express the equations $[X,\Lambda
X]_{\gotico{m}}=0$ as a set of matrix equations.

Denote by $\gotico{m}^{\C}$ the complexification of the tangent
space $\gotico{m}$. We extend $\Lambda$ and the isotropic
representation from $\gotico{m}$ to $\gotico{m}^\C$. Considering
$\gotico{m}^\C_{ij}$, the irreducible submodules of this
representation, we have $\Lambda
\big{|}_{\gotico{m}^\C_{ij}}=\Lambda \big{|}_{\gotico{m}^\C_{ji}}=
\lambda_{ij}\mbox{Id}$.
\\

Denote by $E_{pq}^{ij}$ the $n\times n$  matrix with 1 in position
$(n_1+\ldots + n_{i-1}+p,n_1+\ldots + n_{j-1}+q)$ and zero
elsewhere. The root space associated with the root
$\alpha_{pq}^{ij}:=\varepsilon_p^i-\varepsilon_q^j$ is the complex
span of the matrix $E_{pq}^{ij}$. The matrix subspace
\begin{equation}
M^{ij}=\mbox{ span }\opbra E^{ij}_{pq} \clbra_{0 < p \leq n_i,\, 0<
q \leq n_j}\end{equation} is isomorphic over $\C$ to
$\gotico{m}^\C_{ij}$. Every matrix $A$ can be written as $A=\sum
A^{ij}$, $A^{ij}\in M^{ij}$. With $A^{ij}$ we can associate a matrix
$a_{ij}\in M_{n_i,n_j}(\C);$ specifically, we define
$$A^{ij}=\sum_{p,\,q}z_{pq}\,E^{ij}_{pq}\quad\quad\Rightarrow\quad\quad
a_{ij}=\sum_{p,\,q}z_{pq}E_{pq}.$$ $a_{ij}$ is the only non-trivial
block in $A^{ij}$. Since $A$ is skew-Hermitian we have
$a_{ij}=-a_{ji}^*$.

%More precisaly, $\gotico{m}_{ij}=\bigoplus \gotico{g}_{pq}^{ij}$, where $n_{i-1} \leq p < n_{i}$ e $n_{j-1} \leq q < n_{j}$. Therefore, we write the coeficients of the $\gotico{m}_{ij}$-part of $X$ in the basis $E_{ij}$ how a matrix
%$A_{ij}$ of size $n_i \times n_j$ with complexes entries.

%\begin{obs}
%Given $X\in \gotico{m}\approx T_o \F$, we can write $$X=\sum_{\alpha\in R^+_M}(z_\alpha X_\alpha - \overline{z_\alpha}X_{-\alpha}),$$ where $z_\alpha \in \C$.

%Supose that the coeficients of $\gotico{m}_{ij}$-part of $X$ form the matrix $A_{n_i\times n_j}$. Then the matrix of coeficients of $\gotico{m}_{ji}$-part of $X$ is necessarily $-A^*_{n_j\times n_i}$.
%\end{obs}

\begin{lema}
\label{lemamatrix} Let $i,j,m\in [1,k]$ be distinct. if $X\in
M^{ij}$ and $Y \in M^{jm}$ then $Z=[X,Y]\in M^{im}$. Moreover, if
$X,Y,Z$ are represented by matrix blocks $a\in M_{n_i,n_j}(\C),b\in
M_{n_j,n_m}(\C)$ and $c\in M_{n_i,n_m}(\C)$, repectively, then
$c=ab$.
\end{lema}
\dem This follows from the observation that if $\alpha=\alpha^{i\,
k_1}_{p\, q_1}$ and $\beta=\alpha^{k_2\, j}_{p_2\, q}$ then
$\alpha+\beta$ is a root exactly when $k_1=k_2$ and $q_1=p_2,$ in
which case $\alpha + \beta =\alpha^{i \, j}_{p \, q}$. \fimdem \\

%Writing $X_{\gotico{m}_{ik}}=a_{11}E^{ik}_{11}+\ldots + a_{1\,n_k}E^{ik}_{1\,n_k} + \ldots + a_{n_i\,1}E^{ik}_{n_1\, 1}+ \ldots a_{n_i\,n_k}E^{ik}_{n_1\, n_k}$, \\
%$X_{\gotico{m}_{kj}}=b_{11}E^{kj}_{11}+\ldots + b_{1\,n_j}E^{kj}_{1\,n_j} + \ldots + b_{n_k\,1}E^{kj}_{n_k\, 1}+ \ldots b_{n_k\,n_j}E^{kj}_{n_k\, n_j}$ and using the fact that $[E_\alpha,E_\beta]=m_{\alpha,\beta}\,E_{\alpha+\beta}$ if  $\alpha+\beta$ is root and $0$ otherwise (with $m_{\alpha,\beta}=1$ in our case) we have
%\begin{eqnarray*}
%Z=[X,Y] &=& (a_{11}b_{11}+\ldots + a_{1\,n_k}b_{n_k\,1})E^{ij}_{11} + \ldots + (a_{11}b_{1\,n_k}+\ldots + a_{1\,n_k}b_{n_k\,n_j})E^{ij}_{1\, n_j} \\ &+& (a_{n_i\,1}b_{11}+\ldots + a_{n_i\,n_k}b_{n_k\,1})E^{ij}_{n_i1} + \ldots + (a_{n_i\,1}b_{1\, n_j}+\ldots +a_{n_i\,n_k}b_{n_k\,n_j})E^{ij}_{n_i\,n_j}.
%\end{eqnarray*}
%Therefore, $Z\in \gotico{m}_{ij}$ and the coeficients of $Z$ are the entries of the matrix $A\cdot B$. \fimdem

We can now express the equigeodesic condition in matrix terms.
\begin{teo}
\label{teosistem} Let $X=\sum_{i,j}{X^{ij}}\in \gotico{m}^{\C} $ be
represented by the skew-Hermitian block matrix $A$ with blocks
$a_{ij}\in M_{n_i,n_j}(\C).$ Then $X$ is equigeodesic iff
 \begin{equation}\label{eq.cond}a_{ij}\, a_{jm} = 0\quad\quad\quad
(i,j,m \quad \mbox{distinct}, \quad 1\leq i,j,m \leq
k).\end{equation}
\end{teo}
\dem Let $\Lambda_{ij}$ be the matrix with all-ones in the $ij$ and
$ji$ blocks, and zeros otherwise. Each invariant metric $\Lambda$
has the matrix representation
$\Lambda=\sum{\lambda_{ij}\Lambda_{ij}}$ ($\lambda_{ij} > 0$). It is
clear that the equation $[X,\Lambda X]=0$ ($X\in \gotico{m}$) is
equivalent to $[X,\Lambda_{ij}X]=0$ for all $1\leq i,j \leq k$ ($i
\neq j$). However, a simple calculation based on Lemma
\ref{lemamatrix} shows that the $j$-th block row of
 $[X,\Lambda_{ij}X]=[A,\lambda_{ij}(A^{ij}+A^{ji})]$ consists of the entries
$a_{ji}a_{im}$ $(m\neq i,j)$. Thus, $X$ is equigeodesic iff all
these products vanish.\fimdem

%Let $A$ be the matrix represented by $X$.  The set of all operator $\Lambda$ form a convex cone. Then we have $[X,\Lambda X]=0$ iff $[X,\Lambda_{ij}X]=0$ where $\Lambda_{ij}$ are extremal directions of the cone, i.e., $\Lambda_{ij}$ is a matrix with 1 in positions $(n_1+\ldots + n_{i-1}+p,n_1+\ldots + n_{j-1}+q)$  and $(n_1+\ldots + n_{j-1}+q,n_1+\ldots + n_{i-1}+p)$, $1\leq p < n_i$, $1 \leq q < n_j$ and zero elsewhere. Using the matrix $A$ we have, for each $i,j$ fixed,
%$$[A,\Lambda_{ij} A]=[A,\lambda_{ij}(a_{ij}+a_{ji})]=\left[\sum_{r,\,s}a_{rs} ,\, \lambda_{ij}(a_{ij}+a_{ji})\right]. $$
%The result follow from the linearity of the Lie bracket and from Lemma \ref{lemamatrix}. \fimdem

{~}

According to Theorem \ref{teosistem}, the classification problem for
equigeodesic vectors $X$ in $\F(n;n_1,\\ \cdots,n_s)$ reduces to the
classification problem for the associated skew-Hermitian $n\times n$
matrix $A$ (satisfying the condition $a_{ij}a_{jm}=0),$ up to
conjugation by the unitary subgroup $\hat U:=\oplus_{i=1}^k
U_{n_i}\subset U_n.$ However, as we shall see, closedness of the
associated Killing field really depends on conjugation by the full
unitary group, i.e. depends entirely on the eigenvalues of $A$. We
start with the following definition.

\begin{defi}
We say that a matrix $A$ is {\em essentially diagonal} if $A$ is
permutation-similar to a diagonal matrix, i.e. $A$ contains at most
a single non-zero entry in each row and column.

Analogously, we call $A$ {\em essentially block-diagonal} if $A$
contains at most a single non-zero {\em block entry} $a_{ij}\in
M_{n_i,n_j}(\C)$ in each block-row of size $n_i$ and each column-row
of size $n_j$).
\end{defi}

\noindent In general, neither of these properties implies the other.

\begin{cor}
\label{cormain} $X$ is equigeodesic whenever $A$ is essentially
block-diagonal.
\end{cor}

\noindent Indeed, if $A$ is essentially block-diagonal we have
$a_{ij}a_{jm}=a_{ji}^*a_{jm}=0$ since both $a_{ji}$ and $a_{jm}$
belong to the same block row.\qed

{~}

We remark that each block $a_{ij}$ of $A$ corresponds to one of the
irreducible modules $\gotico{m}_\xi$ defined in the previous
section; moreover, a vector $X$ supported on $\gotico{m}_\xi \oplus
\gotico{m}_\eta$ is essentially block-diagonal exactly when both
$\xi \pm \eta $ are not roots.
\\

\begin{teo}\label{svd}\quad
(i) Every skew-Hermitian matrix $A$ which satisfies (\ref{eq.cond})
is $\hat U$-conjugate to an essentially diagonal matrix $J$.\quad
(ii) The non-zero eigenvalues of $A$ are equal to $\pm i$ times the
absolute value of the non-zero entries of $J$.
\end{teo}

\dem (i) First we discard a few simple cases. The diagonalization of
each block $a_{ij}$ (together with $a_{ji})$ via the SVD algorithm
(singular value decomposition, see e.g. \cite{hj} pp. 157) amounts
to a $\hat U$ conjugation which is non-trivial only in its $i$ and
$j$ block components. If $A$ is essentially block-diagonal, this
step does not change the remaining blocks in $A$, and we may
diagonalize them one by one till an essentially diagonal matrix $J$
is obtained.

If $A$ satisfies (\ref{eq.cond}) but is not essentially
block-diagonal, a slightly more delicate argument is needed. The
skew-symmetry relations $a_{ij}=-a_{ji}^*$ (plus the Fredholm
alternative $Im[A]=Ker[A^*]^\perp$) implies for $i,j,m$ distinct
 $$\begin{array}{cccccc}
 (i) & Im[a_{ji}] & \mbox{and}  & Im[a_{jm}] &
 \mbox{are orthogonal subspaces in} & \C^{n_j},\\
 (ii) & Ker[a_{ij}] & \mbox{and} & Ker[a_{mj}] &
 \mbox{are orthogonal subspaces in} & C^{n_j}.
\end{array}$$
\noindent It follows that all the blocks $a_{ij}$ have orthogonal
cokernels and orthogonal image spaces, hence a single $\hat
U$-conjugation can affect the SVD simultaneously in all of them,
again resulting in an essentially diagonal matrix.

{~}

(ii) $J$ is essentially diagonal and skew-Hermitian, hence it is
permutation-similar to a direct sum of skew-Hermitian $2\times 2$
matrices,
    \begin{equation}\label{ak}
    J_k=\begin{pmatrix} 0&a_k\\ -a_k&0\end{pmatrix}, \hspace{1cm} a_k\geq 0
    \end{equation}
with eigenvalues $\pm i|a_k|.$ Since $A$ and $J$ are similar,
these are also the non-zero eigenvalues of $A$.\fimdem\\

The integers $r_{ij}=rank(a_{ij})$ satisfy the inequalities $\sum_j
r_{ij}\leq n_i.$ These numbers form a partial set of $\hat
U$-conjugation invariance for an equigeodesic vector. A full set of
invariants is supplied by the singular values of each block
$a_{ij}$.

%Each equigeodesic vector $X\in \gotico{m}$ defines uniquely integers $r_{ij}\, (1\leq i<j \leq k)$ via $r_{ij}=\mbox{rank}(a_{ij})$ and we have $\sum_{i\neq j}r_{ij}\leq n_i, \, (i=1,\ldots, k)$.

\begin{exem}
\normalfont (i) Consider the flag manifold $\F(n;n_1,n_2,n_3)$.
According to Theorem \ref{teosistem}, a non-zero vector $X\in
\gotico{m}$, represented by the matrix $A$, is equigeodesic iff the
blocks $a_{12},a_{13},a_{23}$ satisfy $$a_{12}\, a_{23}=0,\quad\quad
a_{13}^* \, a_{12}  =0,\quad\quad a_{23} \, a_{13}^* = 0. $$ $X$ is
essentially block diagonal iff precisely one of these blocks is
non-zero.\\

(ii) Let $X$ be an equigeodesic vector in $\F(n;n-2,1,1).$ If the
corresponding matrix $A$ is not essentially block-diagonal then
$a_{23}=0$ and the vectors $a_{12}$ and $a_{13}$ are non-zero and
orthogonal. Under a simple basis change in $\C^3$ we may assume that
$a_{12}=(a,0,0)^*$ and $a_{13}=(0,b,0)^*.$ Now $A$ is essentially
diagonal, and its non-zero eigenvalues are $\pm ia$ and $\pm ib$.
\\

(iii) We can use the converse process to create complicated
equigeodesic vectors from simple ones. For example, in the flag
manifold $\F(9;3,3,3)$, we start with any essentially diagonal
matrix, say
$$ A= \left(
\begin{tabular}{c|c|c}

\begin{tabular}{ccc}
~0 & ~~~~0 & ~~~0\\ ~0 & ~~~~0 & ~~~0\\ ~0 & ~~~~0 & ~~~0
\end{tabular}
&
\begin{tabular}{ccc}
$\sigma_1$ & 0 & ~0\\ ~0 & $\sigma_2$ & ~0\\ 0 & 0 & ~0
\end{tabular}
&
\begin{tabular}{ccc}
0 & ~0 & 0\\ 0 & ~0 & 0\\ $\sigma_3$ & ~0 & 0
\end{tabular}
\\ \hline

\begin{tabular}{ccc}
$-\sigma_1$ & 0 & 0\\ 0 & $-\sigma_2$ & 0\\ 0 & 0 & 0
\end{tabular}
&
\begin{tabular}{ccc}
~0~ & ~0~ & ~0\\ 0 & 0 & ~0\\ 0 & 0 & ~0
\end{tabular}
&
\begin{tabular}{ccc}
~~~0 & ~~0 & 0\\ ~~~0 & ~~0 & 0\\ ~~~0 & ~~0 & $\sigma_4$
\end{tabular}
\\ \hline

\begin{tabular}{ccc}
~~~0 & ~~~~0 & $-\sigma_3$\\ ~~~0 & ~~~~0 & 0\\ ~~~0 & ~~~~0 & 0
\end{tabular}
&
\begin{tabular}{ccc}
~~~0 & ~~0 & 0\\ ~~~0 & ~~0 & 0\\ ~~~0 & ~~0 & $-\sigma_4$
\end{tabular}
&
\begin{tabular}{ccc}
~~0 & ~~0 & 0\\ ~~0 & ~~0 & 0\\ ~~0 & ~~0 & 0
\end{tabular}

\end{tabular}
\right),$$ where $\sigma_i>0$ for all $i$. Now, each conjugation by
an element of $\hat U:=U\in \mathbb{U}(3) \oplus \mathbb{U}(3)
\oplus \mathbb{U}(3)$ produces a new equigeodesic vector.

%{\bf O próximo exemplo vai entrar na redação final ???? O exemplo ainda está com a notação antiga !!!}

%In the flag manifold $\F(7;3,3,1)$ we have the follow system of equations:
%$$
%\left\{
%\begin{array}{ccc}
%a_{11}\,c_{11}+a_{12}\,c_{21}+a_{13}\,c_{31} &=& 0 \\
%a_{21}\,c_{11}+a_{22}\,c_{21}+a_{23}\,c_{31} &=& 0 \\
%a_{31}\,c_{11}+a_{32}\,c_{21}+a_{33}\,c_{31} &=& 0 \\
%\overline{b_{11}}\,a_{11}+\overline{b_{21}}\,a_{21}+\overline{b_{31}}\,a_{31}&=&0 \\
%\overline{b_{11}}\,a_{12}+\overline{b_{21}}\,a_{22}+\overline{b_{31}}\,a_{32}&=&0 \\
%\overline{b_{11}}\,a_{13}+\overline{b_{21}}\,a_{23}+\overline{b_{31}}\,a_{33}&=&0 \\
%c_{11}\,\overline{b_{11}}+c_{21}\,\overline{b_{21}}+c_{31}\,\overline{b_{31}}&=&0
%\end{array}
%\right. ,
%$$
%This system have 13 ``families'' of solution. An example is:
%$$
%a_{11}= - {\displaystyle \frac {a_{12}\,c_{21} + a_{13}\,c_{31}}{c_{11}}} , \, a_{21}= - {\displaystyle \frac {a_{22}\,c_{21} + a_{23}\,c_{31}}{c_{11}}} ,\, a_{31}= - {\displaystyle \frac {c_{31}\,(a_{12}\,c_{21} + a_{13}\,c_{31})}{c_{11}^{2}}}
%$$
%$$
%a_{32}={\displaystyle \frac {a_{12}\,c_{31}}{c_{11}}} , \,a_{33}={\displaystyle \frac {a_{13}\,c_{31}}{c_{11}}} , \, b_{11}= %- b_{31}\,\overline{({\displaystyle\frac {c_{31}}{c_{11}}} )}, \,b_{21}=0,
%$$
%and the others variables are free.
\end{exem}

In case of the full flag manifold $\F(n)=\F(n;1,\ldots ,1)$, the
blocks $a_{ij}$ are just complex numbers, and $a_{ij}a_{jm}=0$
implies $a_{ij}=0$ or $a_{jm}=0.$ This proves the following result.

\begin{cor}
\label{mainteo}
%The spaces of solutions of the equation (\ref{eqnlegal}) for all invariant metric in $\F(n)$ are:
%\begin{enumerate}
%\item[a)] ${\gotico{u}_{\alpha}}$, $\alpha\in R_M^+$;
%\item[b)] $\gotico{u}_{\alpha_i}\oplus\gotico{u}_{\alpha_j}$ where $\alpha_i,\alpha_j\in R^+_M$, such that $\alpha_i+\alpha_j$ and $\alpha_i-\alpha_j$ are no roots;
%\item[c)] $\gotico{u}_{\alpha_1}\oplus\ldots\oplus\gotico{u}_{\alpha_n}$ where $\alpha_1,\ldots\alpha_n \in R^+_M$ satisfies parwise the condition b).
%\end{enumerate}
$X\in \gotico{m}$ is an equigeodesic vector in $\F(n)$ iff $A$ is
essentially diagonal.
\end{cor}

Thus, for example, the only equigeodesic vectors in $\F(3)$ are
the obvious ones, which belong to the spaces
$\gotico{u}_{12},\gotico{u}_{23},\gotico{u}_{13}$. Observe that
for any two positive roots $\alpha,\beta\in\gotico{sl}(3)$,
necessarily $\alpha+\beta$ or $\alpha-\beta$ is a root.

\section{Closed equigeodesics}

The closeness of a geodesic is a delicate question which involves
global considerations. However, Theorem \ref{svd} allows us to
isolate a set of equigeodesic vectors whose associated homogeneous
equigeodesic is necessarily closed.

First we provide an intuitive description. If $X$ is an equigeodesic
vector, we may assume that its matrix $A=J$ is already in canonical
form, i.e. essentially diagonal. We interpret the permutation
similarity which transforms $J$ into a direct sum of $2\times 2$
matrices as in (\ref{ak}) as an {\em isometric covering} of $\gamma$
by a geodesic $\tilde\gamma$ on a torus; clearly, if the eigenvalues
of $A$ are commensurate, $\tilde\gamma$, hence also $\gamma,$ is
closed. Otherwise, $\tilde\gamma$ is dense on the torus, but
$\gamma$ may or may not be closed on the flag manifold. 

A more rigorous treatment involves not just $\gamma,$ but the whole
Killing field on the flag manifold defined by $\gamma.$ We start
with the following definition.

\begin{defi}
Let $M$ be a manifold. A vector field $T\in \mathfrak{X}(M)$ is closed if
every induced trajectory is closed.
\end{defi}

The following construction of the Killing field $X^*$ associated
with a given vector $X\in \gotico{m}$ is standard (see, for example \cite{arva1}). We define $X^* \in \mathfrak{X}(\F(n;n_1,\ldots,n_k))$ via
$$X^*(pH)=\frac{d}{dt}((\exp tX)\cdot pH)\Big{|}_{t=0}.$$
If $X$ is a homogeneous geodesic vector then the corresponding
homogeneous geodesic $\gamma$ is the trajectory of $X^*$ through the
origin $o$, that is, $\gamma(t)=\phi_t(o)$. If $X^*$ is closed, so
is $\gamma$.

Clearly, $X^*$ is a Killing vector field with respect to any
$SU(n)$-invariant metric. Namely, the generated flow
$\phi_t(\cdot)=L_{(\exp tX)}(\cdot)$ where $L_a$ $(a\in SU(n)$ is
the left translation) is isometric. It follows that the
one-parameter transformation group defined by $\opbra
\phi_t\clbra_{t\in\mathbb{R}}\subset SU(n)$ consists of isometries.
Topologically, this group is either open $(\R)$ or closed $(S^1)$.

%\sout{It is easy to see that $\gamma$ is closed iff $X^*$ is closed
%{\bf seria bom mencionar uma fonte para isto}}.

\begin{teo}(Flores et al, \cite{Piccione})
\label{tpiccione} Let $T$ be a Killing vector field on a Riemannian
manifold $(M,g)$. Then $T$ is closed iff the associated
one-parameter group
 %$\opbra \phi_t\clbra_{t\in\mathbb{R}}\subset Iso(M,g)$
is $S^1$.
\end{teo}

\begin{teo}\label{clos}
Let $X\in \gotico{m}$ be an equigeodesic vector in
$\F(n;n_1,\ldots,n_k)$ represented by the skew-Hermitian matrix $A$.
%\sout{The equigeodesic $\gamma(t)=\exp (tX)
%\cdot o$
Then the corresponding Killing field is closed iff the eigenvalues
of $A$ are commensurate. This in particular implies that the
equigeodesic $\gamma(t)=\exp (tX) \cdot o$ is closed.
\end{teo}

\dem Let $i\theta_1,\ldots,i\theta_n$ be the eigenvalues of $A$. The
1-parameter group of isometries generated by $X^*$ is $\exp
tA=U(\exp tD) U^*$, where $D=diag(i\theta_1,\ldots,i\theta_n)$.
Evidently this group is closed (i.e. diffeomorphic to $S^1$) iff
$\theta_1,\ldots,\theta_n$ are commensurate. On the other hand, by
Theorem \ref{tpiccione} this group is closed iff $X^*$ is
closed.\fimdem

{\em If in Theorem \ref{clos} the eigenvalues are {\em not}
commensurate, the Killing field is not closed, and we do not know
whether $\gamma$ is necessarily open, or dense, in the flag
manifold.}

\begin{obs}
In the case of full flag manifolds, this theorem establishes a connection
 between closed equigeodesics and the equiharmonic non-holomorphic tori
 described by the third author in \cite{neg2}.
\end{obs}

\begin{exem}
In $\F(4)$ consider the equigeodesic vector
$$X=\left(
\begin{array}{cccc}
0 & x & 0 & 0 \\
-\bar{x} & 0 & 0 & 0 \\
0 & 0 & 0  & y \\
0 & 0 & -\bar{y} & 0

\end{array}
\right). $$ The eigenvalues of $X$ are $\pm i|x|,\pm i|y|$. The
equigeodesic determined by $X$ is closed if $x=2$ and $y=3$.
%\sout{ and not closed if $x=1+i$ and $y=1$}.
\end{exem}

In the following simple case we prove that the geodesic, rather than
the Killing field, is closed.

\begin{prop}\label{fnclosed}
In $\F(n),$ every vector of the form $X\in \gotico{u}_{\alpha}$ with
$\alpha\in R_M^+$ is equigeodesic; and the corresponding geodesic,
$\gamma(t)=\exp (tX) \cdot o$, is closed.
\end{prop}
\dem The fact that $X$ is equigeodesic follows from Corollary
\ref{cormain}. Our proof that the equigeodesic is closed is
based on Helgason's proof in \cite{helg} Ch. IV. The subspace
$\gotico{u}_{\alpha}$ is a Lie triple system in the real Lie algebra
$\gotico{su}(n)$; namely, if $X,Y,Z \in \gotico{u}_{\alpha}$ then
$[X,[Y,Z]]\in \gotico{u}_{\alpha}$. Therefore, the subspace
$\gotico{g}^{\prime}=\gotico{u}_{\alpha}+[\gotico{u}_{\alpha},\gotico{u}_{\alpha}]$
is a Lie subalgebra of $\gotico{su}(n)$ which is isomorphic to
$\gotico{su}(2)$. Let $G^{\prime}$ be the connected subgroup of $G$
with Lie algebra $\gotico{g}^{\prime}$ and $M^{\prime}$ the orbit
$G^{\prime}\cdot o$. We can identify $M^{\prime}$ with $G^{\prime}
/(G^{\prime} \cap T),$ a submanifold of $\F(n)$, with
$T_oM^{\prime}=\gotico{u_{\alpha}}$, see \cite{helg} Ch II. Note
that $M^{\prime}=SU(2)/S(U(1)\times U(1))=S^2$ and the induced
Riemannian metric in $M^{\prime}$ is (up to scaling) the normal
metric. This way, geodesics in $\F(n)$ with geodesic vector in
$T_oM^{\prime}$ are of the form $\exp (tX)\cdot o$ where $X\in
u_{\alpha}$, hence are curves in $M^\prime$. Therefore, the
immersion $M^{\prime}\subset\F(n)$ is geodesic at $o$. As
$G^{\prime}$ acts transitively on $M^{\prime},$ it is totally
geodesic in the sense of \cite{helg}. But geodesics in $S^2$ are closed.\fimdem\\

We remark that geodesic curves are 1-dimensional real-harmonic maps,
and in symplectic geometry are closely related to 1-dimensional
complex-harmonic maps. In the proof of Theorem \ref{fnclosed}, an
equigeodesic with tangent vector $X\in\gotico{u}_\alpha$ extends
uniquely to an equiharmonic map $\phi:S^2 \rightarrow \F(n)$ with
tangent space
$\gotico{u}_\alpha.$\\

\section*{Acknowledgments}

We would like to thank Luiz San Martin, Pedro Catuogno and Paolo
Piccione for helpful discussions. This research was partially
supported by FAPESP grant 07/06896-5(Negreiros), CNPq grant
300019/96-3(Cohen) and CAPES/CNPq grant 140431/2009-8(Grama).

\end{document}